\documentclass[12pt,reqno]{amsart}
\usepackage{amssymb}
\usepackage{mathrsfs}
\usepackage[all]{xy}
\newtheorem{theorem}{Theorem}[section]
\newtheorem{lemma}[theorem]{Lemma}
\newtheorem{prop}[theorem]{Proposition}
\newtheorem{corollary}[theorem]{Corollary}
\theoremstyle{definition}
\newtheorem{definition}{Definition}[section]
\theoremstyle{remark}
\newtheorem{remark}{Remark}[section]
\newcommand*{\ptens}[1]{\mathop{\widehat\otimes}_{#1}}
\newcommand*{\Ptens}{\mathop{\widehat\otimes}}
\newcommand*{\lmod}{\mbox{-}\!\mathop{\mathsf{mod}}}
\newcommand*{\rmod}{\mathop{\mathsf{mod}}\!\mbox{-}}
\newcommand*{\bimod}{\mbox{-}\!\mathop{\mathsf{mod}}\!\mbox{-}}
\newcommand*{\lunmod}{\mbox{-}\!\mathop{\mathsf{unmod}}}
\newcommand*{\runmod}{\mathop{\mathsf{unmod}}\!\mbox{-}}
\newcommand*{\biunmod}{\mbox{-}\!\mathop{\mathsf{unmod}}\!\mbox{-}}
\newcommand*{\id}{\mathbf 1}
\newcommand*{\op}{\mathrm{op}}
\newcommand*{\h}{\mathbf h}

\DeclareMathOperator{\Ext}{Ext}
\renewcommand*{\Im}{\mathop{\mathrm{Im}}}
\renewcommand*{\dh}{\mathop{\mathrm{dh}}}
\newcommand*{\injdh}{\mathop{\mathrm{inj.dh}}}
\newcommand*{\wdh}{\mathop{\mathrm{w.dh}}}
\newcommand*{\db}{\mathop{\mathrm{db}}}
\newcommand*{\wdb}{\mathop{\mathrm{w.db}}}
\newcommand*{\ad}{\mathop{\mathrm{ad}}}
\newcommand*{\CC}{\mathbb C}

\newcommand*{\B}{\mathscr B}
\renewcommand*{\H}{\mathscr H}
\newcommand*{\eps}{\varepsilon}
\newcommand*{\xra}{\xrightarrow}
\begin{document}
\title[Approximate characterizations of projectivity
and injectivity]{Approximate characterizations of projectivity
and injectivity for Banach modules}
\author{A. Yu. Pirkovskii}
\thanks{The main results of this paper were obtained whilst
the author held a Royal Society/NATO
Postdoctoral Fellowship at the University of Leeds in 2004--2005.
The author acknowledges with thanks this financial support.
The author is also grateful to H.~G.~Dales and to the School of Mathematics
of the University of Leeds for hospitality.
This work was also partially supported by RFBR grants 05-01-00982 and 05-01-00001.}
\date{}
\begin{abstract}
We characterize projective and injective Banach modules in approximate terms,
generalizing thereby a characterization of contractible Banach algebras
given by F. Ghahramani and R. J. Loy. As a corollary, we show that each
uniformly approximately amenable Banach algebra is amenable.
Some applications to homological dimensions of Banach modules and
algebras are also given.
\end{abstract}
\maketitle

In \cite{GL}, F. Ghahramani and R. J. Loy introduced several
``approximate'' generalizations of amenability and contractibility
for Banach algebras. Among other things, they defined uniformly approximately
contractible and uniformly approximately amenable Banach algebras
and proved that a uniformly approximately contractible algebra is in fact
contractible. In this note we extend their approach to the setting of Banach modules
by showing that projective and injective Banach modules can be
characterized in approximate terms. As a corollary, we
obtain approximate characterizations of biprojective, biflat and
amenable Banach algebras. In particular, we prove
that a uniformly approximately amenable algebra is automatically
amenable, and we obtain an alternative proof of the above-mentioned
result of F.~Ghahramani and R.~J.~Loy concerning uniform approximate contractibility.

\section{Preliminaries}
We begin by recalling some notation and some basic
facts from the homology theory
of Banach algebras. For details, we refer to \cite{X1,X2,Dalesbook}.

Let $A$ be a Banach algebra. We denote by $A\lmod$ (respectively, $\rmod A$,
$A\bimod A$) the category of left Banach $A$-modules
(respectively, right Banach $A$-modules, Banach $A$-bimodules).
For $X,Y$ in any of the above categories, the space of morphisms
from $X$ to $Y$ is denoted by ${_A}\h(X,Y)$ (respectively, $\h_A(X,Y)$,
${_A}\h_A(X,Y)$). The space of bounded linear operators between Banach
spaces $X$ and $Y$ is denoted by $\B(X,Y)$.
If $A$ is unital, then $A\lunmod$, $\runmod A$, and $A\biunmod A$
stand for the corresponding categories of unital Banach
$A$-modules. The unitization of $A$ is denoted by $A_+$.
Recall that $A\bimod A$ is isomorphic to $A^e\lunmod$,
where $A^e=A_+\Ptens A_+^\op$ is the enveloping algebra of $A$.

A chain complex $X_\bullet$ in $A\lmod$ is {\em admissible} if it splits
in the category of Banach spaces.
A morphism $\varkappa$ (respectively, $\sigma$) in $A\lmod$ is an
{\em admissible monomorphism} (respectively, an {\em admissible epimorphism})
if it fits into a short admissible sequence
$0\to X \xra{\varkappa} Y \xra{\sigma} Z\to 0$.
A left Banach $A$-module $P$ is {\em projective} if for each admissible
epimorphism $X\to Y$ in $A\lmod$ the induced map ${_A}\h(P,X)\to {_A}\h(P,Y)$
is onto. Dually,
a left Banach $A$-module $I$ is {\em injective} if for each admissible
monomorphism $Y\to X$ in $A\lmod$ the induced map ${_A}\h(X,I)\to {_A}\h(Y,I)$
is onto.
A left Banach $A$-module $F$ is {\em flat} if for each admissible monomorphism
$Y\to X$ in $\rmod A$ the operator $Y\ptens{A} F\to X\ptens{A} F$ is
topologically injective (i.e., bounded below). Equivalently, $F$ is flat
if and only if the dual module, $F^*$, is injective in $\rmod A$.
The Banach algebra $A$ is {\em biprojective} (respectively, {\em biflat})
if $A$ is projective (respectively, flat) in $A\bimod A$.

An important fact is that the category $A\lmod$ has {\em enough projectives}
and {\em enough injectives}. This means that for each $X\in A\lmod$ there exists
a projective module $P\in A\lmod$ (respectively, an injective module $I\in A\lmod$)
and an admissible epimorphism $P\to X$ (respectively, an admissible monomorphism
$X\to I$). The same is true of the categories $\rmod A$, $A\bimod A$ etc.

A {\em projective resolution} of $X\in A\lmod$ is a chain complex
$P_\bullet=(P_n,d_n)_{n\ge 0}$ in $A\lmod$ consisting of projective modules
together with a morphism $\eps\colon P_0\to X$ such that the augmented
sequence $P_\bullet\xra{\eps} X\to 0$ is an admissible complex.
By reversing arrows, we obtain the definition of an {\em injective resolution}
of $X\in A\lmod$. Since $A\lmod$ has enough projectives and enough injectives,
it follows that each $X\in A\lmod$ has a projective resolution and an injective
resolution.
For $X,Y\in A\lmod$, the space $\Ext^n_A(X,Y)$ is defined to be the $n$th
cohomology of the complex ${_A}\h(P_\bullet, Y)$, where $P_\bullet$
is a projective resolution of $X$.
Equivalently, $\Ext^n_A(X,Y)$ is the $n$th
cohomology of the complex ${_A}\h(X, I^\bullet)$, where $I^\bullet$
is an injective resolution of $Y$.
Note that $\Ext^n_A(X,Y)$ is a seminormed space in a canonical way.
Recall that $X\in A\lmod$ is projective (respectively, injective)
if and only if $\Ext_A^1(X,Y)=0$ (respectively, $\Ext_A^1(Y,X)=0$)
for all $Y\in A\lmod$. For each $X\in A\lmod$ and each $Y\in\rmod A$,
there is a natural topological isomorphism between $\Ext_A^n(X,Y^*)$
and $\Ext_{A^\op}^n(Y,X^*)$.
If $X\in A\bimod A$, then $\Ext^n_{A^e}(A_+,X)$ is topologically isomorphic
to $\H^n(A,X)$, the $n$th
continuous Hochschild cohomology group of $A$ with coefficients in $X$.

We shall use the following explicit description of the space $\Ext_A^1(X,Y)$
(see \cite[7.3.19]{X2}).
Let $A$ be a Banach algebra and let $X,Y\in A\lmod$.
Denote by $Z^1(A\times X,Y)$ the Banach space of all
continuous bilinear maps $f\colon A\times X\to Y$
satisfying
$$
a\cdot f(b,x)-f(ab,x)+f(a,b\cdot x)=0 \qquad (a,b\in A,\; x\in X).
$$
Define $\delta^0\colon \B(X,Y)\to Z^1(A\times X,Y)$ by
$$
(\delta^0 T)(a,x)=a\cdot T(x)-T(a\cdot x)\qquad (a\in A,\; x\in X).
$$
Then we have
\begin{equation}
\label{Ext1}
\Ext^1_A(X,Y)=Z^1(A\times X,Y)/\Im\delta^0.
\end{equation}

A Banach algebra $A$ is {\em contractible} \cite{X2}
(respectively, {\em amenable} \cite{Jhnsn_CBA}) if for each
Banach $A$-bimodule $X$ every continuous derivation from $A$ to $X$
(respectively, to $X^*$) is inner. Equivalently, $A$ is contractible
(respectively, amenable) if and only if $A_+$ is biprojective
(respectively, biflat).

Given Banach spaces $X_1,\ldots ,X_n,Y$, we denote by
$\B^n(X_1\times\cdots\times X_n,Y)$ the Banach space of $n$-linear
continuous maps from $X_1\times\cdots\times X_n$ to $Y$.
The canonical embedding of a Banach space $X$ into $X^{**}$
will be denoted by $i_X$.

\section{Projective and injective Banach modules}

The following lemma is a version of \cite[Lemme 1]{RR} (see also
\cite[0.5.9]{X1}).

\begin{lemma}
\label{lemma:open}
Let $0 \to X^\bullet \to Y^\bullet \to Z^\bullet \to 0$ be
a short exact sequence of cochain complexes of Banach spaces.
Suppose that some map $f$ belonging to the long exact cohomology sequence
$$
\cdots \to H^n(X^\bullet) \to H^n(Y^\bullet) \to H^n(Z^\bullet)
\to H^{n+1}(X^\bullet) \to \cdots
$$
is surjective. Then $f$ is open.
\end{lemma}

For a proof, see \cite[Lemma 1.2]{Pir_Stein}.

\begin{lemma}
\label{lemma:split}
Let
\begin{equation}
\label{XYZ}
0 \to X \xra{\varkappa} Y \xra{\sigma} Z \to 0
\end{equation}
be an admissible sequence in $A\lmod$. Suppose that
the topology on $\Ext^1_A(Z,X)$ is trivial, and that either
$\Ext^1_A(Y,X)=0$ or $\Ext^1_A(Z,Y)=0$. Then \eqref{XYZ} splits
in $A\lmod$.
\end{lemma}
\begin{proof}
First suppose that $\Ext^1_A(Y,X)=0$. Applying the functor ${_A}\h(\,\cdot\, ,X)$
to \eqref{XYZ}, we obtain an exact sequence
\begin{equation*}
0 \to {_A}\h(Z,X) \to  {_A}\h(Y,X) \xra{\varkappa^*} {_A}\h(X,X)
\xra{\delta} \Ext^1_A(Z,X) \to 0.
\end{equation*}
By Lemma~\ref{lemma:open}, $\delta$ is open.
Therefore the triviality of the topology on $\Ext^1_A(Z,X)$ means exactly that
$\varkappa^*$ has dense range. Since the set of invertible elements
in ${_A}\h(X,X)$ is open, there exists $\varphi\in {_A}\h(Y,X)$
such that $\varkappa^*(\varphi)=\varphi\varkappa$ is invertible in ${_A}\h(X,X)$.
Therefore $(\varphi\varkappa)^{-1}\varphi\colon Y\to X$ is a left
inverse of $\varkappa$, and so $\eqref{XYZ}$ splits in $A\lmod$.

Now suppose that $\Ext^1_A(Z,Y)=0$.
The same argument as above applied to the exact sequence
\begin{equation*}
0 \to {_A}\h(Z,X) \to  {_A}\h(Z,Y) \xra{\sigma_*} {_A}\h(Z,Z)
\to \Ext^1_A(Z,X) \to 0
\end{equation*}
yields $\varphi\in {_A}\h(Z,Y)$
such that $\sigma_*(\varphi)=\sigma\varphi$ is invertible in ${_A}\h(Z,Z)$.
Therefore $\varphi(\sigma\varphi)^{-1}\colon Z\to Y$ is a right
inverse of $\sigma$, and so $\eqref{XYZ}$ splits in $A\lmod$.
\end{proof}

\begin{prop}
\label{prop:proj_inj}
Let $A$ be a Banach algebra and let $X\in A\lmod$.

{\upshape (i)} Suppose that the topology on $\Ext^1_A(X,Y)$
is trivial for each $Y\in A\lmod$. Then $X$ is projective.

{\upshape (ii)} Suppose that the topology on $\Ext^1_A(Y,X)$
is trivial for each $Y\in A\lmod$. Then $X$ is injective.
\end{prop}
\begin{proof}
(i) Take an admissible sequence
\begin{equation}
\label{proj_epi}
0 \to Y \to P\to X\to 0
\end{equation}
with $P$ projective. Clearly, \eqref{proj_epi} satisfies the conditions
of Lemma~\ref{lemma:split}. Therefore \eqref{proj_epi} splits in
$A\lmod$, and so $X$ is projective.

(ii) Take an admissible sequence
\begin{equation}
\label{inj_mono}
0 \to X \to I\to Y\to 0
\end{equation}
with $I$ injective. Then the same argument as above shows
that~\eqref{inj_mono} splits in
$A\lmod$, and so $X$ is injective.
\end{proof}

\begin{remark}
\label{rem:unmod}
If $A$ is unital, then there are obvious analogues of
Lemma~\ref{lemma:split} and Proposition~\ref{prop:proj_inj}
for the category $A\lunmod$.
\end{remark}

Let $A$ be a Banach algebra and $X$ a Banach $A$-bimodule.
Given $x\in X$, denote by $\ad_x$ the inner derivation
$A\to X,\; a\mapsto [a,x]$. By definition \cite{GL},
$A$ is {\em uniformly approximately contractible}
(respectively, {\em uniformly approximately amenable})
if for each Banach $A$-bimodule $X$ and each continuous
derivation $D\colon A\to X$ (respectively, $D\colon A\to X^*$)
there exists a net $\{ x_\nu\}$ in $X$ (respectively, in $X^*$)
such that $D=\lim_\nu \ad_{x_\nu}$ in the norm topology.
Equivalently, $A$ is uniformly approximately contractible
(respectively, uniformly approximately amenable)
if for each Banach $A$-bimodule $X$ the topology on
$\H^1(A,X)$ (respectively, on $\H^1(A,X^*)$) is trivial.

Part (i) of the following corollary is due to F.~Ghahramani and R.~J.~Loy~\cite{GL}.

\begin{corollary}
\label{cor:contr_amen}
Let $A$ be a Banach algebra.

{\upshape (i)} Suppose that $A$ is uniformly approximately contractible.
Then $A$ is contractible.

{\upshape (ii)} Suppose that $A$ is uniformly approximately amenable.
Then $A$ is amenable.
\end{corollary}
\begin{proof}
(i) Since $\H^1(A,X)\cong\Ext^1_{A^e}(A_+,X)$, Proposition~\ref{prop:proj_inj} (i)
(cf. also Remark~\ref{rem:unmod})
implies that $A_+$ is projective in $A\bimod A$, i.e., that $A$
is contractible.

(ii) Since $\H^1(A,X^*)\cong\Ext^1_{A^e}(A_+,X^*)\cong\Ext^1_{A^e}(X,A_+^*)$,
Proposition~\ref{prop:proj_inj} (ii) implies that $A_+^*$ is injective in $A\bimod A$,
i.e., that $A_+$ is flat in $A\bimod A$, i.e., that $A$
is amenable.
\end{proof}

\begin{remark}
An alternative proof of part (ii) of Corollary~\ref{cor:contr_amen}
(avoiding injective modules and $\Ext$-spaces)
was recently given in \cite{GLZ}.
\end{remark}

In order to characterize projective and injective Banach modules in approximate
terms, it will be convenient to give the following definition.

\begin{definition}
\label{def:UAM}
Let $A$ be a Banach algebra, and let $X,Y\in A\lmod$.
A {\em uniform approximate morphism} from $X$ to $Y$ is a net $\{ \varphi_\nu\}$
in $\B(X,Y)$ such that $\varphi_\nu(a\cdot x)-a\cdot\varphi_\nu(x)\to 0$
uniformly on bounded subsets of $A$ and $X$.
\end{definition}

Similarly one defines uniform approximate morphisms of right Banach $A$-modules.
If $X,Y\in A\bimod A$, we say that a net $\{\varphi_\nu\}$ in $\B(X,Y)$
is a {\em uniform approximate $A$-bimodule morphism} if it is
a uniform approximate morphism of left and right Banach $A$-modules.

\begin{lemma}
\label{lemma:UAM}
{\upshape (i)} Let $X,Y\in A\bimod A$. A net $\{\varphi_\nu\}$ in $\B(X,Y)$
is a uniform approximate $A$-bimodule morphism if and only if it is
a uniform approximate left $A^e$-module morphism.

{\upshape (ii)} Let $X,Y\in A\lmod$. A net $\{\varphi_\nu\}$
in $\B(X,Y)$ is a uniform
approximate morphism if and only if $\{\varphi_\nu^*\}$ is a uniform
approximate morphism from $Y^*$ to $X^*$.

{\upshape (iii)} If $X,Y\in A\lmod$ and $\{\varphi_\nu\}$ is a uniform
approximate morphism from $X$ to $Y$, then for each pair of morphisms
$\psi\in {_A}\h(Y,Z)$ and $\tau\in {_A}\h(W,X)$ the net $\{\psi\varphi_\nu\tau\}$
is a uniform approximate morphism from $W$ to $Z$.
\end{lemma}
\begin{proof}
(i) The ``if'' part is clear. Conversely, let
$\{\varphi_\nu\}$ be a uniform approximate $A$-bimodule morphism.
For each $\nu$ define a trilinear map $\Phi_\nu\in\B^3(A_+\times X\times A_+,Y)$
by
$$
\Phi_\nu(a,x,b)=\varphi_\nu(a\cdot x\cdot b)-a\cdot\varphi_\nu(x)\cdot b.
$$
Then we have
$$
\| \Phi_\nu(a,x,b) \| \le
\| \varphi_\nu(a\cdot x\cdot b)-a\cdot\varphi_\nu(x\cdot b) \|
+ \| a\| \| \varphi_\nu(x\cdot b)-\varphi_\nu(x)\cdot b \|,
$$
and so $\Phi_\nu\to 0$ in $\B^3(A_+\times X\times A_+,Y)$.
Now it remains to apply the canonical isometric isomorphism
$\B^3(A_+\times X\times A_+,Y)\cong \B^2(A^e\times X,Y)$.

(ii) For each $a\in A$, let $L_a\colon X\to X$ (respectively, $R_a\colon X^*\to X^*$)
denote the left (respectively, right) multiplication by $a$. Since
$R_a=L_a^*$, we have
\begin{equation}
\label{LR}
\| L_a\varphi_\nu-\varphi_\nu L_a\|
= \| \varphi_\nu^* R_a-R_a\varphi_\nu^*\|.
\end{equation}
On the other hand, $\{\varphi_\nu\}$ (respectively, $\{\varphi_\nu^*\}$)
is a uniform approximate morphism if and only if the left-hand side
(respectively, the right-hand side) of \eqref{LR} converges to $0$ uniformly
on bounded subsets of $A$. This proves (ii).

(iii) This is a direct computation.
\end{proof}

\begin{theorem}
\label{thm:proj_appr}
Let $A$ be a Banach algebra. The following properties
of $P\in A\lmod$ are equivalent:

\begin{itemize}
\item[{\upshape (i)}]
$P$ is projective.
\item[{\upshape (ii)}]
The topology on $\Ext^1_A(P,X)$ is trivial for each $X\in A\lmod$.
\item[{\upshape (iii)}]
For each diagram
\begin{equation}
\label{test_proj}
\xymatrix{
& P \ar[d]^\varphi \\
X \ar[r]^\sigma & Y \ar[r] & 0
}
\end{equation}
in $A\lmod$ with $\sigma$ an admissible epimorphism, there
exists a net $\{\varphi_\nu\}$ in ${_A}\h(P,X)$ such that
$\sigma\varphi_\nu\to\varphi$ in the norm topology.
Equivalently, for each admissible
epimorphism $X\to Y$ in $A\lmod$ the induced map ${_A}\h(P,X)\to {_A}\h(P,Y)$
has dense range.
\item[{\upshape (iv)}]
For each diagram \eqref{test_proj} in $A\lmod$ with
$\sigma$ an admissible epimorphism, there
exists a uniform approximate morphism $\{\psi_\nu\}$ in $\B(P,X)$ such that
$\sigma\psi_\nu=\varphi$ for all $\nu$.
\item[{\upshape (v)}]
For each diagram \eqref{test_proj} in $A\lmod$ with
$\sigma$ an admissible epimorphism, there
exists a uniform approximate morphism $\{\psi_\nu\}$ in $\B(P,X)$ such that
$\sigma\psi_\nu\to\varphi$ in the norm topology.
\item[{\upshape (vi)}]
There exists a projective module $F\in A\lmod$, an admissible epimorphism
$\pi\colon F\to P$ and a uniform approximate morphism
$\{\rho_\nu\}$ in $\B(P,F)$ such that
$\pi\rho_\nu\to\id_P$ in the norm topology.
\item[{\upshape (vii)}]
For each $X\in A\lmod$ and each $f\in Z^1(A\times P,X)$, there
exists a net $\{ T_\nu\}$ in $\B(P,X)$ such that
$\delta^0 T_\nu\to f$ in the norm topology.
\end{itemize}
\end{theorem}
\begin{proof}
(i)$\iff$(ii). This follows from Proposition \ref{prop:proj_inj} (i).

(ii)$\iff$(vii). This follows from \eqref{Ext1}.

(i)$\Longrightarrow$(iii). This is clear.

(iii)$\Longrightarrow$(iv). Given diagram \eqref{test_proj},
find a net $\{\varphi_\nu\}$ in ${_A}\h(P,X)$ satisfying (iii).
Let $\rho\in\B(Y,X)$ be a right inverse of $\sigma$.
For each $\nu$, set $\psi_\nu=\varphi_\nu+\rho(\varphi-\sigma\varphi_\nu)$.
Then it is clear that $\sigma\psi_\nu=\varphi$ for all $\nu$,
and
$$
\| \psi_\nu(a\cdot x)-a\cdot\psi_\nu(x)\|
\le  2 \|\rho\| \|\varphi-\sigma\varphi_\nu\| \| a\| \| x\| \to 0
$$
uniformly on bounded subsets of $A$ and $X$.
Therefore the net $\{\psi_\nu\}$ satisfies (iv), as required.

(iv)$\Longrightarrow$(v)$\Longrightarrow$(vi). This is clear.

(vi)$\Longrightarrow$(v). Using the projectivity of $F$, we can complete
each diagram of the form \eqref{test_proj} to a commutative diagram
$$
\xymatrix{
F \ar[r]^\pi \ar@{-->}[d]_\tau & P \ar[d]^\varphi \\
X \ar[r]^\sigma & Y
}
$$
in $A\lmod$. Applying Lemma~\ref{lemma:UAM}, we see that the net
$\psi_\nu=\tau\rho_\nu$ satisfies (v).

(v)$\Longrightarrow$(vii).
Given $f\in Z^1(A\times P,X)$, define $P\times_f X\in A\lmod$
as follows. As a Banach space, $P\times_f X$ coincides with $P\times X$.
For each $a\in A$ and each $(p,x)\in P\times_f X$, set
$$
a\cdot (p,x)=(a\cdot p,a\cdot x-f(a,p)).
$$
Then it is easy to check (cf. \cite[7.2.38]{X2}, \cite[I.1.9]{X1}) that $P\times_f X$
is a left Banach $A$-module.

Clearly, the map
$$
\pi_1\colon P\times_f X\to P,\quad \pi_1(p,x)=p
\qquad (p\in P,\; x\in X)
$$
is an admissible epimorphism. Hence there exists a uniform approximate morphism
$\{\psi_\nu\}$ in $\B(P,P\times_f X)$ such that $\pi_1\psi_\nu\to\id_P$
in the norm topology.
Set $S_\nu=\pi_1\psi_\nu$ and $T_\nu=\pi_2\psi_\nu$,
where $\pi_2\colon P\times_f X\to X$ is the projection onto the second
direct summand. We have
$$
a\cdot\psi_\nu(x)-\psi_\nu(a\cdot x)
=\bigl(a\cdot S_\nu(x),a\cdot T_\nu(x)-f(a,S_\nu(x))\bigr)
-\bigl(S_\nu(a\cdot x),T_\nu(a\cdot x)\bigr).
$$
Therefore
\begin{equation}
\label{dT-f}
(\delta^0 T_\nu-f)(a,x)=\pi_2\bigl(a\cdot\psi_\nu(x)-\psi_\nu(a\cdot x)\bigr)
+f(a,S_\nu(x)-x).
\end{equation}
Since $a\cdot\psi_\nu(x)-\psi_\nu(a\cdot x)\to 0$ uniformly on bounded
subsets of $A$ and $X$ and $S_\nu(x)=(\pi_1\psi_\nu)(x)\to x$ uniformly
on bounded subsets of $X$, it follows from \eqref{dT-f} that
$\delta^0 T_\nu\to f$ in the norm topology, as required.
\end{proof}

\begin{remark}
By replacing in (v) and in Definition~\ref{def:UAM} the norm convergence with the uniform
convergence on compact sets, we obtain the definition of {\em approximately
projective} Banach module introduced by O.~Yu.~Aristov \cite{Ar_appr}.
For approximately projective modules, the equivalences analogous to
$\mathrm{(iv)}\iff\mathrm{(v)}\iff\mathrm{(vi)}\iff\mathrm{(vii)}$ are also valid
(see \cite[6.2]{Ar_appr}).
\end{remark}

\begin{theorem}
\label{thm:inj_appr}
Let $A$ be a Banach algebra. The following properties
of $I\in A\lmod$ are equivalent:

\begin{itemize}
\item[{\upshape (i)}]
$I$ is injective.
\item[{\upshape (ii)}]
The topology on $\Ext^1_A(X,I)$ is trivial for each $X\in A\lmod$.
\item[{\upshape (iii)}]
For each diagram
\begin{equation}
\label{test_inj}
\xymatrix{
0 \ar[r] & Y \ar[r]^\varkappa \ar[d]_\varphi & X\\
& I
}
\end{equation}
in $A\lmod$ with $\varkappa$ an admissible monomorphism, there
exists a net $\{\varphi_\nu\}$ in ${_A}\h(X,I)$ such that
$\varphi_\nu\varkappa\to\varphi$ in the norm topology.
Equivalently, for each admissible
monomorphism $Y\to X$ in $A\lmod$ the induced map ${_A}\h(X,I)\to {_A}\h(Y,I)$
has dense range.
\item[{\upshape (iv)}]
For each diagram \eqref{test_inj} in $A\lmod$ with
$\varkappa$ an admissible monomorphism, there
exists a uniform approximate morphism $\{\psi_\nu\}$ in $\B(X,I)$ such that
$\psi_\nu\varkappa=\varphi$ for all $\nu$.
\item[{\upshape (v)}]
For each diagram \eqref{test_inj} in $A\lmod$ with
$\varkappa$ an admissible monomorphism, there
exists a uniform approximate morphism $\{\psi_\nu\}$ in $\B(X,I)$ such that
$\psi_\nu\varkappa\to\varphi$ in the norm topology.
\item[{\upshape (vi)}]
There exists an injective module $J\in A\lmod$, an admissible monomorphism
$\mu\colon I\to J$ and a uniform approximate morphism
$\{\rho_\nu\}$ in $\B(J,I)$ such that
$\rho_\nu\mu\to\id_I$ in the norm topology.
\item[{\upshape (vii)}]
For each $X\in A\lmod$ and each $f\in Z^1(A\times X,I)$, there
exists a net $\{ T_\nu\}$ in $\B(X,I)$ such that
$\delta^0 T_\nu\to f$ in the norm topology.
\end{itemize}
\end{theorem}

The proof of Theorem~\ref{thm:inj_appr} is similar to that of
Theorem~\ref{thm:proj_appr} and is therefore omitted.

\begin{corollary}
Let $A$ be a Banach algebra. For $X\in A\lmod$, the following
conditions are equivalent:

{\upshape (i)} $X$ is flat;

{\upshape (ii)} for each admissible monomorphism $\varkappa\colon Y\to Z$
in $\rmod A$ the map
$$
(\varkappa\otimes\id_X)^{**}\colon (Y\ptens{A} X)^{**}\to (Z\ptens{A} X)^{**}
$$
is injective.
\end{corollary}
\begin{proof}
$\mathrm{(i)}\Longrightarrow\mathrm{(ii)}$.
If $X$ is flat, then $\varkappa\otimes\id_X$ is topologically injective, and hence
so is $(\varkappa\otimes\id_X)^{**}$.

$\mathrm{(ii)}\Longrightarrow\mathrm{(i)}$. By the adjoint associativity formula
\cite[II.5.21]{X1},
the dual map $(\varkappa\otimes\id_X)^*$ is identified with
$$
\h_A(\varkappa,X^*)\colon\h_A(Z,X^*)\to\h_A(Y,X^*),\quad
\psi\mapsto\psi\varkappa.
$$
Since $(\varkappa\otimes\id_X)^{**}$ is injective, we see that
$\h_A(\varkappa,X^*)$ has dense range. Since this holds for each
admissible monomorphism $\varkappa$, Theorem~\ref{thm:inj_appr} (iii)
implies that $X^*$ is injective in $\rmod A$, i.e., that
$X$ is flat in $A\lmod$.
\end{proof}

\section{Biprojective, biflat and amenable algebras}

Let $A$ be a Banach algebra, and
let $\pi\colon A\Ptens A\to A$ denote the product map.
Recall that $A$ is biprojective (respectively, biflat)
if and only if there exists a morphism $\rho\colon A\to A\Ptens A$
(respectively, $\rho\colon (A\Ptens A)^*\to A^*$) in $A\bimod A$
such that $\pi\rho=\id_A$ (respectively, $\rho\pi^*=\id_{A^*}$);
see \cite[IV.5.6 and VII.2.7]{X1} or \cite[2.8.41]{Dalesbook}.
In this section we give approximate versions of these characterizations.

\begin{corollary}
\label{cor:bipr}
Let $A$ be a Banach algebra. The following conditions are equivalent:

{\upshape (i)} $A$ is biprojective.

{\upshape (ii)} There exists a uniform approximate $A$-bimodule morphism
$\{\rho_\nu\}$ from $A$ to $A\Ptens A$ such that $\pi\rho_\nu\to\id_A$
in the norm topology.
\end{corollary}
\begin{proof}
(i)$\Longrightarrow$(ii). This is clear.

(ii)$\Longrightarrow$(i).
For each $\nu$, denote by
$\rho_\nu^\ell$ the composition of $\rho$ with the canonical embedding of
$A\Ptens A$ into $A_+\Ptens A$. Let $\pi^\ell\colon A_+\Ptens A\to A$ denote the
product map extending $\pi$. Then it is clear that $\pi^\ell$ is an admissible epimorphism,
that $\pi^\ell\rho^\ell_\nu\to\id_A$
in the norm topology, and that $\{\rho^\ell_\nu\}$ is a
uniform approximate morphism from $A$ to $A_+\Ptens A$.
Since $A_+\Ptens A$ is projective in
$A\lmod$, Theorem~\ref{thm:proj_appr} (vi) implies that so is $A$.
Therefore $A\Ptens A_+$ is projective
in $A\bimod A$. Now a similar argument applied to the product map
$\pi^r\colon A\Ptens A_+\to A$ shows that $A$ is projective in
$A\bimod A$.
\end{proof}

\begin{remark}
Weakened forms of condition (ii) were used by Y.~Zhang \cite{Zhang_nilp}
and O.~Yu.~Aristov \cite{Ar_appr} to define various ``approximate'' versions
of biprojectivity.
\end{remark}

\begin{corollary}
\label{cor:bifl}
Let $A$ be a Banach algebra. The following conditions are equivalent:

{\upshape (i)} $A$ is biflat.

{\upshape (ii)} There exists a uniform approximate $A$-bimodule morphism
$\{\rho_\nu\}$ from $(A\Ptens A)^*$ to $A^*$ such that $\rho_\nu\pi^*\to\id_{A^*}$
in the norm topology.

{\upshape (iii)}
There exists a uniform approximate $A$-bimodule morphism
$\{\mu_\nu\}$ from $A$ to $(A\Ptens A)^{**}$
such that $\pi^{**}\mu_\nu\to i_A$ in the norm topology.
\end{corollary}
\begin{proof}
(i)$\iff$(ii). The proof is dual to that of Corollary~\ref{cor:bipr};
cf. \cite[VII.2.7]{X1}

(ii)$\Longrightarrow$(iii). Given $\rho_\nu$, set $\mu_\nu=\rho_\nu^* i_A$.
By Lemma~\ref{lemma:UAM}, $\{\mu_\nu\}$ is a uniform approximate $A$-bimodule morphism.
We have
$$
\pi^{**}\mu_\nu=\pi^{**}\rho_\nu^* i_A=(\rho_\nu\pi^*)^* i_A\to i_A,
$$
as required.

(iii)$\Longrightarrow$(ii).
Given $\mu_\nu$, set $\rho_\nu=\mu_\nu^* i_{(A\Ptens A)^*}$.
By Lemma~\ref{lemma:UAM}, $\{\rho_\nu\}$ is a uniform approximate $A$-bimodule morphism.
We have
$$
\rho_\nu\pi^*=\mu_\nu^* i_{(A\Ptens A)^*}\pi^*
=\mu_\nu^* \pi^{***} i_{A^*}
=(\pi^{**}\mu_\nu)^* i_{A^*}
\to i_A^* i_{A^*} = \id_{A^*},
$$
as required.
\end{proof}

Recall that a Banach algebra $A$ is amenable if and only if there
exists an element $M\in (A_+\Ptens A_+)^{**}$ (called
a {\em virtual diagonal} for $A_+$) such that $a\cdot M=M\cdot a$ for all $a\in A_+$
and $\pi_+^{**}(M)=i_{A_+}(1_+)$, where $\pi_+\colon A_+\Ptens A_+\to A_+$ is the product map
and $1_+$ is the identity of $A_+$
(see \cite[1.3]{Jhnsn_appr} or \cite[VII.2.25]{X1}).
The following corollary shows that a virtual diagonal can be replaced
by a ``uniform approximate virtual diagonal''.

\begin{corollary}
\label{cor:amen}
Let $A$ be a Banach algebra. The following conditions are equivalent:

{\upshape (i)} $A$ is amenable.

{\upshape (ii)} There exists a net $\{ M_\nu\}$ in $(A_+\Ptens A_+)^{**}$
such that $\pi_+^{**}(M_\nu)\to i_{A_+}(1_+)$ and
$a\cdot M_\nu-M_\nu\cdot a\to 0$ uniformly on bounded subsets of $A_+$.
\end{corollary}
\begin{proof}
(i)$\Longrightarrow$(ii). This is clear.

(ii)$\Longrightarrow$(i). Define $\mu_\nu\colon A_+\to (A_+\Ptens A_+)^{**}$
by $\mu_\nu(a)=M_\nu\cdot a$. Obviously, $\mu_\nu$ is a right $A_+$-module morphism.
On the other hand,
$$
\mu_\nu(ab)-a\cdot\mu_\nu(b)=(M_\nu\cdot a-a\cdot M_\nu)\cdot b\to 0
$$
uniformly on bounded subsets of $A_+$, so that $\{\mu_\nu\}$ is a uniform approximate
$A_+$-bimodule morphism. Finally,
$$
\pi_+^{**}\mu_\nu(a)=\pi_+^{**}(M_\nu\cdot a)=\pi_+^{**}(M_\nu)\cdot a\to
i_{A_+}(1_+)\cdot a=i_{A_+}(a)
$$
uniformly on bounded subsets of $A_+$, i.e., $\pi_+^{**}\mu_\nu\to i_{A_+}$
in the norm topology. Now Corollary~\ref{cor:bifl} shows that $A_+$ is biflat,
i.e., that $A$ is amenable.
\end{proof}

\begin{remark}
According to F.~Ghahramani and R.~J.~Loy \cite{GL}, a Banach algebra $A$
is {\em approximately amenable} if for each
Banach $A$-bimodule $X$ and each continuous
derivation $D\colon A\to X^*$ there exists a net $\{ x_\nu\}$ in $X^*$
such that $D=\lim_\nu \ad_{x_\nu}$ in the strong operator topology.
By Theorem~2.1 from \cite{GL}, $A$ is approximately amenable if and only
if it satisfies a condition similar to condition~(ii) from Corollary~\ref{cor:amen}
with uniform convergence replaced by pointwise convergence.
\end{remark}

We end this section with an application to locally compact groups.
Let $G$ be a locally compact group with left Haar measure $\mu$.
If $f$ is a function on $G$, then for each $x\in G$
we denote by ${_x}f$ the function defined by ${_x}f(y)=f(xy)\; (y\in G)$.
We endow $\CC$ with the right Banach $L^1(G)$-module structure
determined by the homomorphism
$$
L^1(G)\to\CC,\quad a\mapsto \int_G a(s) \,d\mu(s).
$$
Using the canonical isomorphism $L^1(G)^*\cong L^\infty(G)$, we
consider $L^\infty(G)$ as a right Banach $L^1(G)$-module. It is easy
to check that for each $f\in L^\infty(G)$ and each $a\in L^1(G)$ we have
$f\cdot a=\tilde a*f$, where $\tilde a(x)=\Delta(x^{-1}) a(x^{-1})$
and $\Delta$ is the modular function on $G$.

\begin{corollary}
Let $G$ be a locally compact group. Then $G$ is amenable if and only if
there exists a net $\{ m_\nu\}$ in $L^\infty(G)^*$ such that
$m_\nu(1)\to 1$ and $m_\nu({_x}f)-m_\nu(f)\to 0$ uniformly on $G$
and on bounded subsets of $L^\infty(G)$.
\end{corollary}
\begin{proof}
The ``only if'' part is clear. Conversely, let $\{ m_\nu\}$ be a net
with the indicated properties. We claim that $\{ m_\nu\}$ is a
uniform approximate morphism from $L^\infty(G)$ to $\CC$.
Indeed, for each $f\in L^\infty(G)$ and each $a\in L^1(G)$ we have
\begin{align*}
| m_\nu(f\cdot a)-m_\nu(f)\cdot a|
&=\left| m_\nu\left(\int_G \tilde a(x) ({_{x^{-1}}}f-f)\, d\mu(x)\right)\right|\\
&=\left| \int_G \tilde a(x) m_\nu({_{x^{-1}}}f-f)\, d\mu(x)\right|\\
&\le \| a\| \sup_{x\in G} |m_\nu({_{x^{-1}}}f-f)| \to 0
\end{align*}
uniformly on bounded subsets of $L^\infty(G)$ and $L^1(G)$.
Therefore $\{ m_\nu\}$ is a
uniform approximate morphism from $L^\infty(G)$ to $\CC$.

Define a linear map
$i\colon\CC\to L^\infty(G)$ by $i(1)=1$. Clearly, $i$ is an admissible
monomorphism in $\rmod L^1(G)$, and $m_\nu i\to \id_\CC$.
Since $L^1(G)$ is projective in
$L^1(G)\lmod$ (see \cite[IV.2.17]{X1}), it follows that $L^\infty(G)$ is injective
in $\rmod L^1(G)$. Now Theorem~\ref{thm:inj_appr} (vi) implies that
$\CC$ is injective in $\rmod L^1(G)$. By \cite[VII.2.33]{X1}, this is equivalent
to the amenability of $G$.
\end{proof}

\section{Homological dimensions}

In this final section, we generalize Proposition \ref{prop:proj_inj}
and Corollary \ref{cor:contr_amen} to higher $\Ext$-groups.

Let $A$ be a Banach algebra.
Recall that the {\em projective homological dimension}
(respectively, the {\em injective homological dimension}) of
$X\in A\lmod$ is the least integer $n$ with the property
that $\Ext^{n+1}_A(X,Y)=0$ (respectively, $\Ext^{n+1}_A(Y,X)=0$)
for each $Y\in A\lmod$. The {\em weak homological dimension}
of $X\in A\lmod$ is the least integer $n$ with the property
that $\Ext^{n+1}_A(X,Y^*)=0$ for each $Y\in \rmod A$.
The projective (respectively, injective, flat)
homological dimension of $X$ is denoted by $\dh_A X$
(respectively, $\injdh_A X$, $\wdh_A X$).
Obviously, $\wdh_A X=\injdh_{A^\op} X^*$ for each $X\in A\lmod$.
Note that $X$ is projective
(respectively, injective, flat) if and only if
$\dh_A X=0$ (respectively, $\injdh_A X=0$, $\wdh_A X=0$).

For a Banach algebra $A$, the numbers $\db A=\dh_{A^e} A_+$
and $\wdb A=\wdh_{A^e} A_+$ are called
{\em the homological bidimension} and the
{\em weak homological bidimension} of $A$, respectively.
Note that $A$ is contractible (respectively, amenable) if and only
if $\db A=0$ (respectively, $\wdb A=0$).

\begin{prop}
Let $A$ be a Banach algebra and let $X\in A\lmod$.

{\upshape (i)} Suppose that the topology on
$\Ext^{n+1}_A(X,Y)$ is trivial for each $Y\in A\lmod$. Then $\dh_A X\le n$.

{\upshape (ii)} Suppose that the topology on
$\Ext^{n+1}_A(Y,X)$ is trivial for each $Y\in A\lmod$. Then $\injdh_A X\le n$.

{\upshape (iii)} Suppose that the topology on
$\Ext^{n+1}_A(X,Y^*)$ is trivial for each $Y\in \rmod A$. Then $\wdh_A X\le n$.
\end{prop}
\begin{proof}
(i) Take an admissible sequence
$$
0 \to Z \to P_{n-1} \to \cdots \to P_0 \to X\to 0
$$
with $P_0,\ldots ,P_{n-1}$ projective. A standard argument (cf. \cite[II.5.4]{X1})
shows that $\Ext^{n+1}_A(X,Y)\cong\Ext^1_A(Z,Y)$ for each $Y\in A\lmod$,
the isomorphism being topological by Lemma~\ref{lemma:open}.
By Proposition~\ref{prop:proj_inj} (i), $Z$ is projective.
Therefore $\Ext^{n+1}_A(X,Y)\cong\Ext^1_A(Z,Y)=0$ for each $Y\in A\lmod$,
i.e., $\dh_A X\le n$, as required.

(ii) Take an admissible sequence
$$
0 \to X \to I_0 \to \cdots \to I_{n-1} \to Z\to 0
$$
with $I_0,\ldots ,I_{n-1}$ injective.
Then the same argument as above shows that $Z$ is injective.
Therefore $\Ext^{n+1}_A(Y,X)\cong\Ext^1_A(Y,Z)=0$ for each $Y\in A\lmod$,
i.e., $\injdh_A X\le n$, as required.

(iii) This follows from (ii) applied to $X^*\in A^\op\lmod$.
\end{proof}

\begin{corollary}
Let $A$ be a Banach algebra.

{\upshape (i)} Suppose that the topology on $\H^{n+1}(A,X)$
is trivial for each $X\in A\bimod A$. Then $\db A\le n$.

{\upshape (ii)} Suppose that the topology on $\H^{n+1}(A,X^*)$
is trivial for each $X\in A\bimod A$. Then $\wdb A\le n$.
\end{corollary}

\noindent\textbf{Acknowledgments. }The author is grateful to
H.~G.~Dales, F.~Ghahramani and A.~Ya.~Helemskii for valuable discussions.

\vspace*{20mm}
\begin{flushleft}
\scshape\small
Department of Nonlinear Analysis and Optimization\\
Faculty of Science\\
Peoples' Friendship University of Russia\\
Mikluho-Maklaya 6\\
117198 Moscow\\
RUSSIA

\medskip
{\itshape Address for correspondence:}\\

\medskip\upshape
Krupskoi 8--3--89\\
Moscow 119311\\
Russia

\medskip
{\itshape E-mail:} {\ttfamily pirkosha@sci.pfu.edu.ru, pirkosha@online.ru}
\end{flushleft}

\begin{thebibliography}{88}
\bibitem{Ar_appr}
Aristov, O. Yu.
{\em On the approximation of flat Banach modules by free Banach modules}.
Math. Sbornik \textbf{196} (2005), no.~11 (in Russian).
\bibitem{Dalesbook}
Dales, H. G.
{\em Banach algebras and automatic continuity}.
London Mathematical Society Monographs, New Series, 24.
Oxford Science Publications,
The Clarendon Press, Oxford University Press,
New York, 2000.
\bibitem{GL}
Ghahramani, F., Loy, R. J.
{\em Generalized notions of amenability}.
J. Funct. Anal. \textbf{208} (2004). no.~1, 229--260.
\bibitem{GLZ}
Ghahramani, F., Loy, R. J., Zhang, Y.
{\em Generalized notions of amenability} II. In preparation.
\bibitem{X1}
Helemskii, A. Ya. {\itshape The Homology of Banach and Topological Algebras},
Moscow University Press, 1986 (Russian); English transl.: Kluwer Academic
Publishers, Dordrecht, 1989.
\bibitem{X2}
Helemskii, A. Ya.
{\itshape Banach and Polynormed Algebras: General Theory,
Representations, Homology},
Nauka, Moscow, 1989 (Russian);
English transl.: Oxford University Press, 1993.
\bibitem{Jhnsn_CBA}
Johnson, B. E.
{\em Cohomology in Banach algebras},
Mem. Amer. Math. Soc. \textbf{127} (1972).
\bibitem{Jhnsn_appr}
Johnson, B. E.
{\em Approximate diagonals and cohomology of certain annihilator
Banach algebras}.
Amer. J. Math. \textbf{94} (1972). no.~3, 685--698.
\bibitem{Pir_Stein}
Pirkovskii, A. Yu.
{\em On certain homological properties of Stein algebras}.
J. Math. Sci. (New York) \textbf{95} (1999), no.~6, 2690--2702.
\bibitem{RR}
Ramis, J.-P., Ruget, G.
{\em Complexe dualisant et th\'eor\`emes de dualit\'e en g\'eom\'etrie
analytique complexe}.
Publ. Math. I.H.E.S. \textbf{38} (1970), no.~2, 77--91.
\bibitem{Zhang_nilp}
Zhang, Y.
{\em Nilpotent ideals in a class of Banach algebras}.
Proc. Amer. Math. Soc. \textbf{127} (1999), no.~11, 3237--3242.
\end{thebibliography}
\end{document}